\documentclass[letterpaper,11pt]{article}
\pdfoutput=1

\usepackage[totalwidth=480pt, totalheight=680pt]{geometry}
\usepackage{setspace}

\usepackage{amsmath, array}
\usepackage{mathtools}
\usepackage{fancyvrb}
\usepackage{dsfont}
\usepackage{rsfso}

\usepackage[english,greek]{babel}
\usepackage{arydshln}

\usepackage[usenames,dvipsnames]{color}
\usepackage{color}
\usepackage{pstricks,framed}
\usepackage{amssymb}
\usepackage{mathtools}
\usepackage{tikz}
\usetikzlibrary{chains}

\usepackage{amssymb}
\usepackage{amsthm}
\usepackage{amsfonts}

\usepackage{chngcntr}
\usepackage[retainorgcmds]{IEEEtrantools}

\usepackage[colorlinks=true, linktoc=page, citecolor=blue, urlcolor=blue]{hyperref}
\usepackage{mathdots}

\VerbatimFootnotes

\begin{document}
\selectlanguage{english}
\title{{\color{Blue}\textbf{New Minimal Hypersurfaces in $\mathbb{R}^{(k+1)(2k+1)}$ and $\mathbb{S}^{2k^2+3k}$}}\\[12pt]} \date{}
\author{\textbf{Jens Hoppe},$^1$ \textbf{Georgios Linardopoulos},$^{2,3}$ \textbf{{O. Teoman} Turgut}$^4$\footnote{E-mails:\href{mailto:hoppe@math.kth.se}{hoppe@kth.se}, \href{mailto:glinard@inp.demokritos.gr}{glinard@inp.demokritos.gr}, \href{mailto:turgutte@boun.edu.tr}{turgutte@boun.edu.tr}.}\\[12pt]
$^1$ Department of Mathematics, Royal Institute of Technology,\\
Lindstedtsv\"{a}gen 25, S-10044 Stockholm, Sweden. \\[6pt]
$^2$ Institute of Nuclear and Particle Physics, N.C.S.R.\ ``Demokritos",\\
153 10 Agia Paraskevi, Greece.\\[6pt]
$^3$  Department of Physics, National and Kapodistrian University of Athens,\\
Zografou Campus, 157 84 Athens, Greece. \\[6pt]
$^4$ Department of Physics, Bo\u{g}azi\c{c}i University, \\
34342, Bebek Istanbul, Turkey.}
\maketitle
\begin{abstract}
\normalsize{\noindent We find a class of minimal hypersurfaces $\mathcal{H}_k$ as the zero level set of Pfaffians, resp.\ determinants of real $2k+2$ dimensional antisymmetric matrices. While $\mathcal{H}_1$ and $\mathcal{H}_2$ are congruent to the quadratic cone $x_1^2 + x_2^2 + x_3^2 - x_4^2 - x_5^2 - x_6^2 = 0$ resp.\ Hsiang's cubic $\mathfrak{su}\left(4\right)$ invariant in $\mathbb{R}^{15}$, $\mathcal{H}_{k>2}$ (special harmonic $SO\left(2k+2\right)$-invariant cones of degree $\geq 4$) seem to be new.\\}
\end{abstract}
\noindent Minimal hypersurfaces in $\mathbb{R}^n$ that can be defined by an algebraic equation have a long history, and some classification results are known in degree 2 and 3 (see e.g.\ \cite{Hsiang1967a}, \cite{Tkachev10d}). In this short note we add to the few known (classes of) explicit examples defined as zero level sets a class of hypersurfaces $\mathcal{H}_k$ ($k = 1,2,3,\ldots$), proven to be minimal in $\mathbb{R}^{\left(k+1\right)\left(2k+1\right)}$ by showing that the expression
\begin{IEEEeqnarray}{c}
\boldsymbol{\nabla}\left(\frac{\boldsymbol{\nabla} u}{\left|\nabla u\right|}\right) = \frac{\triangle u}{\left|\nabla u\right|} - \frac{\boldsymbol{\nabla} u}{\left(\nabla u\right)^2}\boldsymbol{\nabla}\left(\left|\nabla u\right|\right), \label{MinimalSurfaceEquation}
\end{IEEEeqnarray}
which is known to be proportional to the mean curvature (of a hypersurface described as a level set by $u\left(x\right) = 0$), vanishes on $u\left(x\right) = 0$ (cp.\ equations \eqref{DefinitionP}--\eqref{Result1} below). \\[6pt]
\indent For any $m = 2k+2 = 2\ell$ dimensional real antisymmetric matrix $M$ consider the Pfaffian:
\begin{IEEEeqnarray}{c}
p\left(M\right):= \frac{1}{2^{\ell}\ell!}\,\epsilon_{a_1 a_2 \ldots a_{2\ell}}M_{a_1a_2}\ldots M_{a_{2\ell-1}a_{2\ell}} \label{PfaffianDefinition}
\end{IEEEeqnarray}
with
%
\begin{IEEEeqnarray}{c}
M = \Big\{M_{ab}\Big\} = \left(\begin{array}{ccccc} 0 & x_1 & \ldots & x_{2\ell - 2} & x_{2\ell - 1} \\[6pt] \vdots && \ddots & \vdots &\vdots \\[6pt] -x_{2\ell - 2} && \ldots & 0 & x_{\ell\left(2\ell-1\right)} \\[6pt]  -x_{2\ell - 1} && \ldots & -x_{\ell\left(2\ell-1\right)} & 0 \end{array}\right) = -M^{T}.
\end{IEEEeqnarray}
$p\left(M\right)$ is a homogeneous polynomial $P_{\ell}\left(\textbf{x}\right)$ of degree $\ell$ in the $n:= \ell\left(2\ell - 1\right) = \left(k+1\right)\left(2k+1\right)$ real variables $x_1, x_2, \ldots x_n$. \\[6pt]
\indent Define
\begin{IEEEeqnarray}{c}
\mathcal{H}_k := \left\{\textbf{x} \in \mathbb{R}^n \ \Big| \ P_{k+1}\left(\textbf{x}\right) = 0\right\}. \label{DefinitionP}
\end{IEEEeqnarray}
While it is obvious that $P_{\ell}\left(\textbf{x}\right) =: u\left(\textbf{x}\right)$ is harmonic (each $x_i$ appears only linearly), hence the first term on the rhs of \eqref{MinimalSurfaceEquation} being trivially = 0, it is also not difficult to prove that
\begin{IEEEeqnarray}{c}
u_i u_j u_{ij} = \rho\left(\textbf{x}\right) u, \label{Property1}
\end{IEEEeqnarray}
with $\rho$ non singular (in particular on $u\left(x\right) = 0$), where for simplicity we denote $\partial u/\partial x_i$ by $u_i$ and $\partial^2 u/\partial x_i\partial x_j$ by $u_{ij}$. \\[6pt]
%
\indent Hence $\mathcal{H}_k$ is a minimal hypersurface, as $\left(\nabla u\right)^2 \triangle u - u_i u_j u_{ij}$ (which is proportional to the mean curvature, cp.\ e.g.\ \cite{Hsiang1967a}) vanishes on $u = 0$. To prove \eqref{Property1} note that because of $p\left(\textbf{x}\right)^2 = \det M$,
\begin{IEEEeqnarray}{c}
\frac{\partial p}{\partial M_{ab}} = p \, M^{ba}
\end{IEEEeqnarray}
($M^{ba}$ denoting the $ba$ matrix element of $M^{-1}$) from which it easily follows that
\begin{IEEEeqnarray}{c}
S_{ab,cd}:= \frac{\partial^2 p}{\partial M_{ab}\partial M_{cd}} = p \left(M^{ba}M^{dc} - M^{bc}M^{da} + M^{bd}M^{ca}\right), \label{DefinitionS}
\end{IEEEeqnarray}
as for any invertible matrix
\begin{IEEEeqnarray}{c}
\frac{\partial M^{ba}}{\partial M_{cd}} = - M^{be} \frac{\partial M_{ef}}{\partial M_{cd}} M^{fa}.
\end{IEEEeqnarray}
While it is obvious (cp.\ \eqref{PfaffianDefinition}) that \eqref{DefinitionS} is polynomial (of degree $\ell-2$) in the matrix elements $M_{..}$, the crucial observation is that it is totally antisymmetric in $\left[abcd\right]$, so that when calculating the l.h.s.\ of \eqref{Property1}, $\partial p/\partial M_{ab} \cdot \partial p/\partial M_{cd}$ can be replaced by its antisymmetric part (obtained by fully anti-symmetrizing $p^2 M^{ba}M^{dc}$)
\begin{IEEEeqnarray}{c}
\frac{1}{3} p^2 \left(M^{ba}M^{dc} - M^{bc}M^{da} + M^{bd}M^{ca}\right) = \frac{1}{3} p \, S_{ab,cd}\,; \label{Antisymmetrizer1}
\end{IEEEeqnarray}
hence
\begin{IEEEeqnarray}{c}
u_i u_j u_{ij} = \frac{1}{12} p \, \text{Tr}\left[S^2\right] = \rho\left(\textbf{x}\right)u, \label{Result1}
\end{IEEEeqnarray}
with $\rho\left(\textbf{x}\right) = \text{Tr}\left[S^2\right]/12$ and Tr$\left[S^2\right] \equiv S_{ab,cd}\,S_{cd,ab}$. \\[6pt]
\indent While the method of this proof works equally well for determinants of matrices with unconstrained entries\footnote{Note \cite{Tkachev10b} where a (different) proof was given for the unconstrained determinental family.} (the resulting index structure $M^{ba}M^{dc} - M^{bc}M^{da}$ still has enough [anti]symmetry to conclude a result similar to \eqref{Result1}), one should perhaps explicitly note that the crucial symmetry argument leading to \eqref{Antisymmetrizer1}--\eqref{Result1} would \textit{not} work for (determinants $\Delta$ of) symmetric matrices, as for them one would only get
\begin{IEEEeqnarray}{c}
\frac{\partial \Delta}{\partial M_{ab}} = 2 \Delta M^{ba},
\end{IEEEeqnarray}
while
\begin{IEEEeqnarray}{c}
\frac{\partial^2 \Delta}{\partial M_{ab}\partial M_{cd}} = 2\Delta \left(M^{ba}M^{dc} - M^{bc}\frac{M_{ef}}{\partial M_{cd}}M^{fa}\right) = 2\Delta \left(2M^{ba}M^{dc} - M^{bc}M^{da} - M^{bd}M^{ca}\right)
\end{IEEEeqnarray}
has no particular (anti)symmetry w.r.t.\ $a \leftrightarrow c$ or $b \leftrightarrow d$. Therefore the two first derivative factors of the product $4\Delta^2 M^{ba} M^{dc}$ can in that case \textit{not} be written as $\Delta$ times a polynomial in the original matrix variables (one should of course also note that the determinant $\Delta$ will no longer be harmonic). \\[6pt]
\indent Noting that $P_2\left(\textbf{x}\right) = x_1 x_6 - x_2 x_5 + x_3 x_4$ (which upon $45^\circ$ rotations in the $\left(16\right)$, $\left(25\right)$ and $\left(34\right)$ planes is easily seen to be a standard minimal quadratic cone) let us now discuss the case $\ell = 3$ ($m = 6$, resp.\ $k = 2$): \\
\begin{IEEEeqnarray}{c}
M = \left(\begin{array}{cccccc} 0 & x_1 & x_2 & x_3 & x_4 & x_5 \\[6pt] -x_1 & 0 & x_6 & x_7 & x_8 & x_9 \\[6pt]  -x_2 & -x_6 & 0 & x_{10} & x_{11} & x_{12} \\[6pt] -x_3 & -x_7 & -x_{10} & 0 & x_{13} & x_{14} \\[6pt] -x_4 & -x_8 & -x_{11} & -x_{13} & 0 & x_{15} \\[6pt] -x_5 & -x_9 & -x_{12} & -x_{14} & -x_{15} & 0 \end{array}\right) \label{Antisymmetric6x6}
\end{IEEEeqnarray} \\[6pt]
for which
\begin{IEEEeqnarray}{ll}
p_3\left(M\right) = \frac{1}{48} \epsilon_{abcdef} M_{ab} M_{cd} M_{ef} & = x_1\left(x_{10}x_{15} - x_{11}x_{14} + x_{12}x_{13}\right) - x_2\left(x_{7}x_{15} - x_{8}x_{14} + x_{9}x_{13}\right) + \nonumber \\[6pt]
& + x_3\left(x_{6}x_{15} - x_{8}x_{12} + x_{9}x_{11}\right) - x_4\left(x_{6}x_{14} - x_{7}x_{12} + x_{9}x_{10}\right) + \nonumber \\[6pt]
& + x_5\left(x_{6}x_{13} - x_{7}x_{11} + x_{8}x_{10}\right) = P_3\left(x_1, \ldots, x_{15}\right),
\end{IEEEeqnarray}
\begin{IEEEeqnarray}{ll}
\frac{\partial p_3}{\partial M_{ab}} = \frac{1}{8} \epsilon^{abcdef} M_{cd} M_{ef}
\end{IEEEeqnarray}
and
\begin{IEEEeqnarray}{ll}
\frac{\partial^2 p_3}{\partial M_{ab}\partial M_{cd}} = \frac{1}{2} \epsilon^{abcdef} M_{ef}.
\end{IEEEeqnarray}
While for cubic cones some kind of classification does exist (due to Tr$\left[S^2\right] = 24\textbf{x}^2$, Tr$\left[S^3\right] = 48 P_3\left(\textbf{x}\right)$, $\mathcal{H}_2$ would be called an exceptional eigencubic of type $\left(0,8\right)$ in the notation of \cite{Tkachev10d}; see also \cite{NadirashviliTkachevVladut14}) and abstract arguments (isomorphism between $\mathfrak{so}\left(6\right)$ and $\mathfrak{su}\left(4\right)$, and each of the cubics invariant under the respective groups) necessitate a congruence of $\mathcal{H}_2$ with Hsiang's hypersurface in $\mathbb{R}^{15}$ (which is obtained by setting to zero the trace of the third power of a general anti-hermitian $4\times4$ matrix), it is reassuring to verify the correspondence explicitly: using that the Dynkin diagrams of
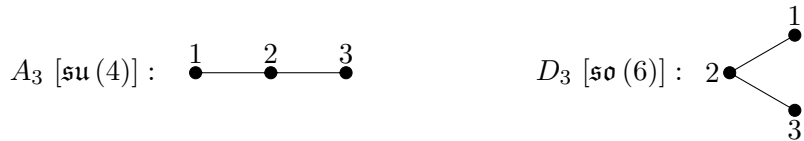
\begin{center}
\begin{tikzpicture}
\node[] at (0,0) {$A_3 \ [\mathfrak{su}\left(4\right)]:$};
\draw[fill=black]
(1.5,0)
      circle [radius=.08] node [above] {$1$} --
(2.5,0)
      circle [radius=.08] node [above] {$2$} --
      node [midway,above] {}
(3.5,0)
      circle [radius=.08] node [above] {$3$}
      node [midway,above] {}

;\node[] at (7,0) {$D_3 \ [\mathfrak{so}\left(6\right)]:$};
\draw[fill=black]

(8.6,0) circle [radius=.08]

(8.6,0) node [left] {$2$}

(8.6,0) --++ (30:1)
      circle [radius=.08] node [above] {$1$}

(8.6,0) --++ (-30:1)
      circle [radius=.08] node [below] {$3$}
;

\end{tikzpicture}
\end{center}
are identical one gets, by taking in the latter case the $6\times6$ matrices $E_{23}-E_{65}$, $E_{12}-E_{54}$ and $E_{26}-E_{35}$ as corresponding to the positive simple roots of \\
\begin{IEEEeqnarray}{c}
L:= \left\{X\in \mathfrak{gl}\left(6,\mathbb{C}\right) \; \Big| \; X^T = -MXM\right\}, \quad M = \left(\begin{array}{cc} 0 & \mathds{1} \\ \mathds{1} & 0 \end{array}\right)
\end{IEEEeqnarray} \\
whereas the usual $E_{12}'$, $E_{23}'$ and $E_{34}'$ for $\mathfrak{su}\left(4\right)$, and using the explicit map  \\
\begin{IEEEeqnarray}{c}
X \rightarrow Y:= \frac{1}{2}\left(\begin{array}{cc} \mathds{1} & 0 \\ 0 & -i\cdot\mathds{1} \end{array}\right)\left(\begin{array}{cc} \mathds{1} & \mathds{1} \\ -\mathds{1} & \mathds{1} \end{array}\right)\ X \ \left(\begin{array}{cc} \mathds{1} & -\mathds{1} \\ \mathds{1} & \mathds{1} \end{array}\right)\left(\begin{array}{cc} \mathds{1} & 0 \\ 0 & i\cdot\mathds{1} \end{array}\right)
\end{IEEEeqnarray} \\
for the correspondence of $L$ with the [complexification of] $\mathfrak{so}\left(6\right)$ (as the Lie algebra of real antisymmetric $6\times6$ matrices), one finds
\begin{IEEEeqnarray}{l}
Z= \frac{i}{2}\left(\begin{array}{cccc} x_3+x_8-x_{12} & x_9+x_{11} & -x_5-x_{10} & x_7-x_4 \\ x_9+x_{11} & x_3+x_{12}-x_8 & x_4+x_7 & x_{10}-x_5 \\ -x_5-x_{10} & x_4+x_7 & x_8+x_{12}-x_3 & x_{11}-x_9 \\ x_7-x_4 & x_{10}-x_5 & x_{11}-x_9 & -x_3-x_8-x_{12} \end{array}\right) + \\[6pt] \nonumber
\hspace{6cm} + \frac{1}{2}\left(\begin{array}{cccc} 0 & x_{15}+x_{6} & -x_{14}-x_{2} & x_1-x_{13} \\ -x_{15}-x_6 & 0 & x_1+x_{13} & x_2-x_{14} \\ x_{14}+x_2 & -x_1-x_{13} & 0 & x_6-x_{15} \\ -x_1+x_{13} & -x_2+x_{14} & -x_6+x_{15} & 0  \end{array}\right) \qquad
\end{IEEEeqnarray} \\
as being the (anti-hermitian) element in $\mathfrak{su}\left(4\right)$ corresponding to \eqref{Antisymmetric6x6}. It is then straightforward to explicitly verify that
\begin{IEEEeqnarray}{c}
\frac{i}{3}\cdot\text{Tr}\left[Z^3\right] = p_3\left(M\right).
\end{IEEEeqnarray}
This in particular gives a nice geometric understanding of Hsiang's minimal cone in $\mathbb{R}^{15}$, as the vanishing of the Pfaffian is equivalent to the columns (rows) of equation \eqref{Antisymmetric6x6} being linearly dependent. It also points out that Hsiang's minimal hypersurface is singular not only at the origin, but also when \emph{four} of the eigenvalues of \eqref{Antisymmetric6x6}, $\left(\pm \lambda_1, \pm \lambda_2, \pm \lambda_3\right)$ vanish; that singular submanifold is given by the vanishing of the coefficient of $\lambda^2$ in the characteristic equation of \eqref{Antisymmetric6x6} (a quartic), namely

\small\begin{IEEEeqnarray}{ll}
x_{3}^2 x_{6}^2 + x_{4}^2 x_{6}^2 + x_{5}^2 x_{6}^2 + x_{2}^2 x_{7}^2 + x_{4}^2 x_{7}^2 + x_{5}^2 x_{7}^2 + x_{2}^2 x_{8}^2 + x_{3}^2 x_{8}^2 + x_{5}^2 x_{8}^2 + x_{2}^2 x_{9}^2 + x_{3}^2 x_{9}^2 + x_{4}^2 x_{9}^2 + x_{1}^2 x_{10}^2 + x_{4}^2 x_{10}^2 + \nonumber \\[6pt]
x_{5}^2 x_{10}^2 + x_{8}^2 x_{10}^2 + x_{9}^2 x_{10}^2 + x_{1}^2 x_{11}^2 + x_{3}^2 x_{11}^2 + x_{5}^2 x_{11}^2 + x_{7}^2 x_{11}^2 + x_{9}^2 x_{11}^2 + x_{1}^2 x_{12}^2 + x_{3}^2 x_{12}^2 + x_{4}^2 x_{12}^2 + x_{7}^2 x_{12}^2 + x_{8}^2 x_{12}^2 + \nonumber \\[6pt]
x_{1}^2 x_{13}^2 + x_{2}^2 x_{13}^2 + x_{5}^2 x_{13}^2 + x_{6}^2 x_{13}^2 + x_{9}^2 x_{13}^2 + x_{12}^2 x_{13}^2+ x_{1}^2 x_{14}^2 + x_{2}^2 x_{14}^2 + x_{4}^2 x_{14}^2 + x_{6}^2 x_{14}^2 + x_{8}^2 x_{14}^2 + x_{11}^2 x_{14}^2 + x_{1}^2 x_{15}^2 + \nonumber \\[6pt]
x_{2}^2 x_{15}^2 + x_{3}^2 x_{15}^2 + x_{6}^2 x_{15}^2 + x_{7}^2 x_{15}^2 + x_{10}^2 x_{15}^2 - 2 x_{2} x_{3} x_{6} x_{7} - 2 x_{2} x_{4} x_{6} x_{8} - 2 x_{3} x_{4} x_{7} x_{8} - 2 x_{2} x_{5} x_{6} x_{9} - 2 x_{3} x_{5} x_{7} x_{9} - \nonumber \\[6pt]
2 x_{4} x_{5} x_{8} x_{9} + 2 x_{1} x_{3} x_{6} x_{10} - 2 x_{1} x_{2} x_{7} x_{10} - 2 x_{1} x_{2} x_{8} x_{11} - 2 x_{3} x_{4} x_{10} x_{11} - 2 x_{7} x_{8} x_{10} x_{11} - 2 x_{1} x_{2} x_{9} x_{12} - \nonumber \\[6pt]
2 x_{3} x_{5} x_{10} x_{12} - 2 x_{7} x_{9} x_{10} x_{12} - 2 x_{4} x_{5} x_{11} x_{12} - 2 x_{8} x_{9} x_{11} x_{12} - 2 x_{1} x_{3} x_{8} x_{13} - 2 x_{2} x_{3} x_{11} x_{13} - 2 x_{6} x_{7} x_{11} x_{13} - \nonumber \\[6pt]
2 x_{1} x_{3} x_{9} x_{14} - 2 x_{2} x_{3} x_{12} x_{14} - 2 x_{6} x_{7} x_{12} x_{14} - 2 x_{4} x_{5} x_{13} x_{14} - 2 x_{8} x_{9} x_{13} x_{14} - 2 x_{11} x_{12} x_{13} x_{14} - 2 x_{1} x_{4} x_{9} x_{15} - \nonumber \\[6pt]
2 x_{2} x_{4} x_{12} x_{15} - 2 x_{6} x_{8} x_{12} x_{15} - 2 x_{3} x_{4} x_{14} x_{15} - 2 x_{7} x_{8} x_{14} x_{15} - 2 x_{10} x_{11} x_{14} x_{15} + 2 x_{1} x_{4} x_{6} x_{11} + 2 x_{1} x_{5} x_{6} x_{12} + \nonumber \\[6pt]
2 x_{1} x_{4} x_{7} x_{13} + 2 x_{2} x_{4} x_{10} x_{13} + 2 x_{6} x_{8} x_{10} x_{13} + 2 x_{1} x_{5} x_{7} x_{14} + 2 x_{2} x_{5} x_{10} x_{14} + 2 x_{6} x_{9} x_{10} x_{14} + 2 x_{1} x_{5} x_{8} x_{15} + \nonumber \\[6pt]
2 x_{2} x_{5} x_{11} x_{15} + 2 x_{6} x_{9} x_{11} x_{15} + 2 x_{3} x_{5} x_{13} x_{15} + 2 x_{7} x_{9} x_{13} x_{15} + 2 x_{10} x_{12} x_{13} x_{15} = 0, \label{SingularSubset}
\end{IEEEeqnarray} \normalsize \\
a non-empty, lower dimensional subset of (the 14-dimensional hypersurface) $\mathcal{H}_2$. To see why the singular set \eqref{SingularSubset} is at most 9-dimensional, note that the equation $\nabla P_3\left(x_1,\ldots, x_{15}\right) = 0$ can be solved for at least 6 out of the 15 variables $x_i$. Therefore the regular subset of $P_3 = 0$\footnote{That is the set $\left\{\textbf{x} \in \mathbb{R}^{15}: P_{3}\left(\textbf{x}\right) = 0 \wedge \nabla P_{3}\left(\textbf{x}\right) \neq 0\right\}$.} is 14-dimensional. The argument trivially goes through for all the higher-dimensional hypersurfaces $P_{k+1} = 0$.
\section*{Acknowledgements}
We would like to thank Jaigyoung Choe for initiating our interest in minimal cones---and Per B\"{a}ck, Aksel Bergfeldt, as well as Martin
Bordemann, for discussions. G.L. and T.T. would like to thank the Royal Institute of Technology for hospitality and support during the collaboration, while J.H. thanks Bogazici University and the Istanbul Center for Mathematical Sciences.
\bibliographystyle{halpha}
\bibliography{Bibliography}
\end{document}